\documentclass[12pts]{amsart}
\usepackage{amsmath}
\usepackage{amsfonts}
\usepackage{flafter}
\usepackage{amscd,amsthm}
\usepackage{amsfonts}
\usepackage{graphicx}
\usepackage{verbatim}
\usepackage{amssymb,amsthm,amsfonts,amstext}
\usepackage{amsmath}
\usepackage{graphicx}
\usepackage{calc}
\usepackage{color}

\usepackage{xspace}

\usepackage{psfrag}
\usepackage{epsfig}

\usepackage{a4wide}

\makeatletter \@addtoreset{equation}{section} \makeatother

\renewcommand\thetable{\thesection.\@arabic\c@table}
\newtheorem{theorem}{Theorem}[section]

\newtheorem{corollary}[theorem]{Corollary}

\newcommand{\cqd}{\hfill$\sqcup\!\!\!\!\sqcap\bigskip$}

\newcommand\ASEP{\textsc{asep}\xspace}
\def\TASEP{\textsc{tasep}\xspace}

\def\ZZ{\mathbb{Z}}
\def\NN{\mathbb{N}}
\def\cR{\mathcal{R}}
\def\cB{\mathcal{B}}

\newcommand{\un}[1]{\underline #1}

\newcommand{\teta}{\tilde\eta}

\usepackage{wasysym}
\newcommand{\parti}{{\textcolor[gray]{0.2}{\CIRCLE}}}
\newcommand{\second}{{\circledast}}
\newcommand{\hole}{{\Circle}}

\newcommand{\hathole}{{\hat\hole}}
\newcommand{\hatparti}{{\hat\parti}}

\newcommand{\starpair}{$^*$pair\xspace}
\newcommand{\starpairs}{$^*$pairs\xspace}
\newcommand{\starhole}{$^*$hole\xspace}
\newcommand{\starparticle}{$^*$particle\xspace}

\newcommand\coalesce{{collide }}

\begin{document}

\title[Collision probabilities in the ASEP]{Collision probabilities in the
  rarefaction fan\\ of asymmetric exclusion processes} \date{August 19, 2008}
\author{Pablo A. Ferrari, Patr{\i}cia Gon\c calves and James B. Martin}

\begin{abstract}
  We consider the one-dimensional asymmetric simple exclusion process (\ASEP) in
  which particles jump to the right at rate $p\in(1/2,1]$ and to the left at
  rate $1-p$, interacting by exclusion. In the initial state there is a finite
  region such that to the left of this region all sites are occupied and to the
  right of it all sites are empty. Under this initial state, the hydrodynamical
  limit of the process converges to the rarefaction fan of the associated
  Burgers equation. In particular suppose that the initial state has first-class
  particles to the left of the origin, second-class particles at sites 0 and 1,
  and holes to the right of site 1. We show that the probability that the two
  second-class particles eventually collide is $(1+p)/3p$, where a
  \emph{collision} occurs when one of the particles attempts to jump over the
  other.  This also corresponds to the probability that two \ASEP processes,
  started from appropriate initial states and coupled using the so-called
  ``basic coupling'', eventually reach the same state.  We give various other
  results about the behaviour of second-class particles in the \ASEP. In the
  totally asymmetric case ($p=1$) we explain a further representation in terms
  of a multi-type particle system, and also use the collision result to derive
  the probability of coexistence of both clusters in a two-type version of the
  corner growth model.
\end{abstract}
\renewcommand{\subjclassname}{\textup{2000} Mathematics Subject Classification}

\address{\noindent IME, Instituto de Matem\'atica e Estat\'{\i}stica da Universidade de S\~ao
  Paulo, Rua do Mat\~ao 1010, Cidade Universit\'aria, S\~ao Paulo, Brasil
\newline$ $
\newline
$\mbox{ }\,\,\,\,$ Department of Statistics, University of Oxford, 1 South Parks Road, Oxford OX1 3TG, UK
\newline$ $}
\email{pablo@ime.usp.br, patg@ime.usp.br, martin@stats.ox.ac.uk}
\maketitle
\section{Introduction}


In the one-dimensional \textit{asymmetric simple exclusion process} (\ASEP{}),
particles perform continuous-time random walks on $\mathbb{Z}$ with rate
$p\in(1/2,1]$ of jumping to the right and $q=1-p$ of jumping to the left.
Particles interact by exclusion; attempted jumps to occupied sites are suppressed.
The resulting process $\eta_{t}$ is Markov in the state space
$\{\hole, \parti\}^{\mathbb{Z}}$, where for a site $x\in{\mathbb{Z}}$, $\eta_{t}(x)=\parti$
indicates that the site $x$ is occupied at time $t$, and $\eta_t(x)=\hole$ indicates that 
it is empty. If $p=1$ then all jumps are to the right and the process is known as the
\emph{totally asymmetric simple exclusion process} (\TASEP{}).



A convenient way to construct the process is using a family of independent Poisson 
processes $\big\{N(x,x\pm1), x\in\ZZ\big\}$, where for each $x$, 
$N(x,x+1)$ has rate $p$ and $N(x,x-1)$ has rate $q$. At the time of 
a point in the process $N(x,y)$, if there is a particle at site $x$, it 
attempts to jump to site $y$; the jump succeeds if site $y$ is empty.
Two versions of the process started from different initial conditions can be jointly constructed under the 
\textit{basic coupling} by using the same $N(x,y)$ for both versions. 

A \emph{second-class particle}, denoted $\second$, is a particle that
interacts with particles like a hole, and interacts with holes
like a particle. That is, a second-class particle jumps
right at rate $p$ if there is a hole on its right, and jumps
left at rate $q$ if there is a hole on its left;
if the site on its left contains a particle, they exchange positions
at rate $p$, and if the site on its right contains a particle,
they exchange positions at rate $q$.

A second-class particle can be regarded as a discrepancy between two systems evolving together
under the basic coupling. If the systems start from two initial conditions which coincide except at a single site,
where one has a particle and the other has a hole, then there will be exactly one discrepancy between
the two systems for all time, and this discrepancy moves as a second-class particle.

Consider an initial state in which every negative site is occupied by
a particle, every positive site has a hole, and there is a single second-class particle
at the origin. Let $X(t)$ be the position of the second-class particle at time $t$.
For the \TASEP{}, 
Ferrari and Kipnis \cite{F.K.} proved
that $X(t)/t$ converges in distribution as $t\rightarrow{\infty}$ to a uniform
random variable on the interval $[-1,1]$. The extension of the method of \cite{F.K.} to 
the \ASEP{} is straightforward, but for completeness we give the result and its proof
(indeed it can be extended much further, to more general asymmetric
exclusion processes and to a larger class of initial distributions).
For the particular case of the \TASEP{}, almost sure convergence has been proved
by Mountford and Guiol~\cite{M.G.}, Ferrari and Pimentel~\cite{F.P.}, and
Ferrari, Martin and Pimentel~\cite{F.M.P.}.

We then consider systems with two second-class particles. We will be
particularly interested in the event that one of the second-class particles
attempts a jump onto the site of the other one. Various conventions are possible
for what happens when such a jump is attempted; the particles may annihilate, or
coalesce, or one may have priority over the other, or no change at all may
occur. (See also Section \ref{growthproofsec} for discussion
of another relevant possibility). 
The difference between these conventions will not be important for our
purposes, since we concentrate on the probability that such a jump
attempt happens at all. For convenience we will refer to such an event as the
\emph{collision} of the two second-class particles.

We start the process from a deterministic configuration, denoted
${\hat\parti}{\un\second}{\second}{\hat\hole}$,
that has all negative sites occupied
by (first-class) particles, the origin and site 1 occupied by second-class particles,
and all sites to the right of site 1 empty 
(in this notation the origin is underlined and $\hat\parti$ and $\hat\hole$ denote the semi-infinite
sequences $\dots\parti\parti\parti$ and $\hole\hole\hole\dots$ respectively).
We show that the probability that the two second-class particles \coalesce is $\frac{1+p}{3p}$.
If instead the particles start at distance two apart, from the configuration
${\hat\parti}{\un\second}{\hole}{\second}{\hat\hole}$, then the corresponding probability is $\frac{1+2p^2}{6p^2}$.

These probabilities can also be interpreted as ``coupling probabilities''
for systems with only one type of particle. If two one-type {\ASEP}s are 
started from the initial conditions ${\hat\parti}{\un{\parti}}{\hole}{\hat\hole}$ and 
$\hat\parti{\un\hole}{\parti}{\hat\hole}$
and evolve together under the basic coupling,
then the probability that they eventually reach the same state is $\frac{1+p}{3p}$.
The second probability $\frac{1+2p^2}{6p^2}$ applies to the case where the initial conditions
of the two systems are ${\hat\parti}{\un\parti}{\hole}{\hole}{\hat\hole}$ and 
${\hat\parti}{\un\hole}{\hole}{\parti}{\hat\hole}$.

In the particular case of the \TASEP{}, where jumps only occur in one direction, 
various interesting alternative representations are possible. 
Consider a multi-type system in which the initial condition has a particle labelled 
$i$ at site $i$, for all $i\in\ZZ$. Lower-numbered particles
have priority over higher-numbered particles; that is, the 
particles start out in reverse order of priority and every jump consists of a particle
labelled $i$ exchanging places with a particle labelled $j$ to its right, for some $i<j$. 
Whenever the particle at site $x$ has priority over the particle at site $x+1$,
such an exchange between sites $x$ and $x+1$ occurs at rate 1. 
From the coupling/collision probabilities described above, one can deduce
that the probability that at some point particle $i$ overtakes particle $i+1$ is $2/3$,
while the probability that particle $i$ overtakes both particle $i+1$ and particle $i+2$ is $1/2$.
 
The \TASEP{} with two second-class particles can also be mapped onto a 
two-type version of the ``corner growth model'' associated to directed last-passage percolation.
This model has also recently been studied by Coupier and Heinrich \cite{CouHei}.
The event that the two particles stay separated for ever in the \TASEP{} corresponds
to the event that both clusters in the growth process grow unboundedly, without 
one surrounding the other. In \cite{CouHei} it is proved that this event has positive probability. 
Using the correspondence with the \TASEP{}, we obtain that the probability is in fact 1/3; 
this gives perhaps the first example of a model where such a ``coexistence probability'' 
can be precisely calculated. 

In Section \ref{s2} we define notation 
and state our main results for the general $\ASEP$, along with some
discussion. In Section \ref{TASEPsec} we discuss the particular case of the $\TASEP$ and
explain the representation as a multi-type particle system, and in 
Section \ref{growthsec}
we describe the two-type competition growth model.
Section \ref{s3} collects together some known
results on hydrodynamics and couplings. Proofs of our results are given in Sections
\ref{ss4}, \ref{ss5} and \ref{ss6}.

\section{Statement of Results}
\label{s2}
Let $p\in(1/2,1]$ and write $q=1-p$.
We define the two-type \ASEP{} as follows.
For $i,j\in\{\hole,\second,\parti\}$,
we say that $i$ has priority over $j$ and write
$i\prec j$ if $i=\parti, j\in\{\hole,\second\}$ or if $i=\second, j=\hole$. Let $\Xi=\{\hole,\second,\parti\}^{\mathbb Z}$ 
and
consider the generator given on cylinder functions
$f:\Xi\rightarrow{\mathbb{R}}$ by:
\begin{equation}
  Lf(\eta)=\sum_{x\in{\mathbb{Z}}}
\Big\{
p{\mathbf 1} _{\{\eta(x)\prec\eta(x+1)\}}
\,[f(\eta^{x,x+1})-f(\eta)]
+q{\mathbf 1} _{\{\eta(x)\prec\eta(x-1)\}}]
\,[f(\eta^{x,x-1})-f(\eta)]
\Big\}, \label{a30gen}
\end{equation}
where
\[ \eta^{x,y}(z)=\left\{
\begin{array}{rl} \label{etaexchanges}
\eta(z), & \mbox{if $z\neq{x,y}$}\\
\eta(y), & \mbox{if $z=x$}\\
\eta(x), & \mbox{if $z=y$}
\end{array}.
\right.
\]



If the initial configuration $\eta$
belongs to $\{\hole,\parti\}^{\mathbb Z}$, we recover the \ASEP{}
with a single type of particle.
If $p=1$ the attempted jumps are only to the right, and we have the \TASEP{}.

The existence of this process can be proven using Liggett's construction \cite{L.}
or using the Harris graphical construction as in \cite{FM1,FM2}. See also the comments
in Section \ref{s3} about the basic coupling.

Our results concern systems in which all sites to the left of a certain point
contain first-class particles, and all sites to the right of a certain point
contain holes.  (This property is preserved by the dynamics). As
  mentioned in the introduction, we use a compact notation for such states:
$\hat\parti=\dots{\parti}{\parti}{\parti}$ denotes a semi-infinite string of
first-class particles, $\hat\hole={\hole}{\hole}{\hole}\dots$ denotes a
semi-infinite string of holes, and an underline indicates the position of the
origin. For example, ${\hat{\parti}}{\un\second}{\second}{\hat\hole}$ denotes
the configuration defined by
\[
({\hat{\parti}}{\un\second}{\second}{\hat\hole})(x):= \left\{
  \begin{array}{ll}
    \parti,&\hbox{if } x<0;\\
    \second,&\hbox{if } x\in\{0,1\};\\
    \hole,&\hbox{if } x>1\,.
\end{array}
\right.
\]
(In fact the absolute position of the origin is rarely important for the
results, and sometimes we may omit the underline, but it is often convenient to fix it).


The following result extends to the \ASEP{} a result proved for the \TASEP{}
by Ferrari and Kipnis \cite{F.K.}:
\begin{theorem} \label{2cptheorem}
Consider the \ASEP{} with $p\in(\frac12, 1]$, starting from the state ${\hat\parti}{\un\second}{\hat\hole}$.
Let $X(t)$ be the position of the second-class particle at time $t$.
Then $X(t)/{t}\to\mathcal{U}_p$
in distribution as $t\to\infty$, where $\mathcal{U}_p$
is uniformly distributed on the interval $[-(p-q), p-q]$.
\end{theorem}
We note that in fact the same method of proof can be used to prove a much more
general result. Already in \cite{F.K.}, the initial condition can be taken to be
the product measure with density $\lambda$ on negative sites and density $\rho$
on positive sites, with $\lambda>\rho$.  In addition one can
extend to more general asymmetric exclusion processes, in which the particles
(still interacting via exclusion) carry out more general continuous-time random
walks with some drift $\gamma>0$.  The convergence is then to a uniform
distribution on the interval $[\gamma(1-2\lambda), \gamma(1-2\rho)]$.

In the particular case of the \TASEP{}, the convergence has been shown to hold
almost surely \cite{M.G.}, \cite{F.P.}, \cite{F.M.P.}.


We now turn to systems with a pair of second-class
particles. The following theorem again extends to the \ASEP{} results proved for the \TASEP{} in \cite{F.K.}.
The proofs are analogous; we include them for completeness.

\begin{theorem} \label{t2} Consider the \ASEP{} with $p\in(\frac12,1]$,
starting from $\hat\parti{\un\second}{\second}{\hat\hole}$.
Let $\tau$ be the first time that a jump is attempted by one of the 
second-class particles onto the site of the other (which may be infinity).
For $t<\tau$, let $X(t)$ and $Y(t)$ be the positions of the two second-class
particles at time $t$, with $X(t)<Y(t)$.
\begin{itemize}
\item[(a)]
For any $r\in[-(p-q), (p-q)]$,
\begin{equation}
\label{separationprob}
\lim_{t\to\infty}P\big(\tau>t, X(t)\leq rt<Y(t)) = \frac{(p-q)^2-r^2}{4p(p-q)}.
\end{equation}
\item[(b)]
\begin{equation} \label{eq1}
\lim_{t\rightarrow{\infty}}\frac{1}{t}E\Big[\big(Y(t)-X(t)\big)I(\tau>t)\Big]=\frac{(p-q)^2}{3p}.
\end{equation}
\end{itemize}
\end{theorem}
In display (3.7) of \cite{F.K.}, the limit \eqref{separationprob} was stated
for the case $p=1$, but it has the extra term ``$-ru(r,1)$'' which does not
belong there. 
 
The following is our main result:
\begin{theorem} \label{t1} 
Consider an \ASEP{} with 
$p\in(\frac12,1]$, containing two second-class particles. 
Let $\tau$ be the first time that a jump is attempted by one of the second-class
particles onto the site of the other (which may be infinity).
\begin{itemize}
\item[(a)]
If the initial condition is ${\hat\parti}{\un\second}{\second}{\hat\hole}$, then 
\begin{equation}
\label{mainprob}
P(\tau<\infty)=\frac{1+p}{3p}.
\end{equation}
\item[(b)]
If instead the initial condition is ${\hat\parti}{\un\second}{\hole}{\second}{\hat\hole}$, then
\begin{equation}
P(\tau<\infty)=\frac{1+2p^2}{6p^2}.
\end{equation}
\end{itemize}
\end{theorem}
In the \TASEP{} ($p=1$) these probabilities are $2/3$ in part (a) (an upper
bound of $3/4$ was obtained by \cite{F.K.} in this case) and $1/2$ in part (b).
The value in part (b) is obtained from that in part (a) using a short argument
based on conditioning on the first jump and some symmetries of the process under
this initial configuration.

These probabilities also represents a coupling probability in the context of the
\ASEP{} with only one class of particle.  Consider two processes started from
the states ${\hat\parti}{\un\parti}{\hole}{\hat\hole}$ and
${\hat\parti}{\un\hole}{\parti}{\hat\hole}$, and coupled using the \textit{basic
  coupling}, under which jumps are attempted at the same sites at the same times
in both processes (see Section \ref{s3} for details). At time 0 there are
discrepancies between the two processes at sites 0 and 1; these discrepancies
move around, each behaving as a second-class particles, and at
time $\tau$, the first time that one jumps onto the site of the other, 
the two processes reach the same state
and stay coupled thereafter. Hence the probability that they ever reach the same
state is $\frac{1+p}{3p}$, or $2/3$ for the case of the \TASEP.  If instead the
initial conditions are ${\hat\parti}{\un\parti}{\hole}{\hole}{\hat\hole}$ and
${\hat\parti}{\un\hole}{\hole}{\parti}{\hat\hole}$ then the two discrepancies
start at sites 0 and 2, and the coupling probability is $\frac{1+2p^2}{6p^2}$,
giving $1/2$ for the case of the \TASEP.

Our proof of Theorem \ref{t1} involves estimating the expected number of
particles which are to the right of the second-class particle at time $t$
in a process started from the initial condition $\hat\parti{\un{\second}}{\hat\hole}$.
As a corollary of these methods, we can deduce the following result,
which relates to the convergence of the environment around a second-class particle:
\begin{corollary}
\label{2cpcorollary}
Let $\eta_t, t\geq 0$ be the \ASEP{} with $p\in(\frac12,1]$ started from the state
${\hat\parti}{\un\second}{\hat\hole}$. Let $X(t)$ be the position of the second-class particle at
time $t$. Then 
\begin{align}
\label{correq1}
\lim_{t\to\infty}E \eta_t\big(X(t)-1\big)
&=\frac12-\frac{p-q}6,
\\
\intertext{and}
\label{correq2}
\lim_{t\to\infty}E \eta_t\big(X(t)+1\big)
&=\frac12+\frac{p-q}6.
\end{align}
\end{corollary}






\section{Multi-class \TASEP{}}
\label{TASEPsec}
When $p=1$, the system is ``totally asymmetric'' and various interesting connections can be drawn.

We consider a system with infinitely many classes of particle. Initially, site $i$
contains a particle labelled $i$, for all $i\in\ZZ$. A particle with label $i$ 
has priority over a particle with label $j$ if $i<j$. 
If the particle at site $x$ has priority over the particle at site $x+1$, 
the two exchange places with rate 1. Since the process is totally 
asymmetric, once two particles have exchanged places, 
they can never swap back again, so the higher-priority particle
stays ahead of the lower-priority particle thereafter.

Let $X_i(t)$ be the position at time $t$ of the particle with label $i$ 
(so $X_i(0)=i$ for all $i$). Consider $i=0$ for example.
Particle $0$ has lower priority than
all the particles starting to its left, and higher priority than all 
the particles starting to its right. Hence $X_0(t)$ behaves 
just as the position $X(t)$ of the second-class particle starting 
from the configuration ${\hat\parti}{\un\second}{\hat\hole}$. The same is true for any other $i$
(shifted by a constant distance $i$ to the right).
In particular, from the strong law of large numbers mentioned after Theeorem \ref{2cptheorem},
we know that the limits 
\begin{equation}\label{Ddef}
D_i=\lim X_i(t)/t
\end{equation}
exist almost surely for all $i$,
and all have the uniform distribution on $[-1,1]$.

Now consider the particles labelled $0$ and $1$. Again they see stronger particles
to their left and weaker particles to their right, so their
behaviour corresponds to that of the pair of second-class particles
starting from the configuration ${\hat\parti}{\un\second}{\second}{\hat\hole}$. 
In particular, we know from Theorem \ref{t1}(a) that the probability 
of a jump attempt between these two particles is $2/3$. 
Since the process is totally asymmetric, the only way for this 
to happen is that particle $0$ attempts to jump right 
onto particle $1$ at some point. In this case the
two particles exchange places. The same applies for any pair of particles $i$ and $i+1$.
So we obtain the following result giving the probability that a particle
eventually overtakes its original right-hand neighbour:
\begin{corollary}
For any $i$,
$P\big(X_i(t)>X_{i+1}(t)\textup{ for some } t\big)
=2/3$.
\end{corollary}

From Theorem \ref{t2}(a), we can obtain further
information about the joint distibution of the ``asymptotic speeds'' $D_i$ and $D_{i+1}$
of neighbouring particles $i$ and $i+1$, as defined at (\ref{Ddef}). We get that
\[
P\big(D_i<r, D_{i+1}>r\big)=\frac{(1-r)(1+r)}4
=P(D_i<r)P(D_{i+1}>r)
\]
for all $r\in[-1,1]$.
This expression would be consistent with $D_i$ and $D_{i+1}$ being independent;
however, since we know that particle $i$ overtakes particle $i+1$ 
with probability $2/3$, such independence does not in fact hold.

From the second part of Theorem \ref{t1} we can make the following
deduction about the probability that a particle overtakes
both of the two particles starting just to its right:
\begin{corollary}\label{distancetwocorr}
For any $i$,
$P\big(X_i(t)>\max\{X_{i+1}(t), X_{i+2}(t)\}\textup{ for some } t\big)
=1/2$.
\end{corollary}

\begin{figure}[h]
\centering
\includegraphics[width=0.5\linewidth]{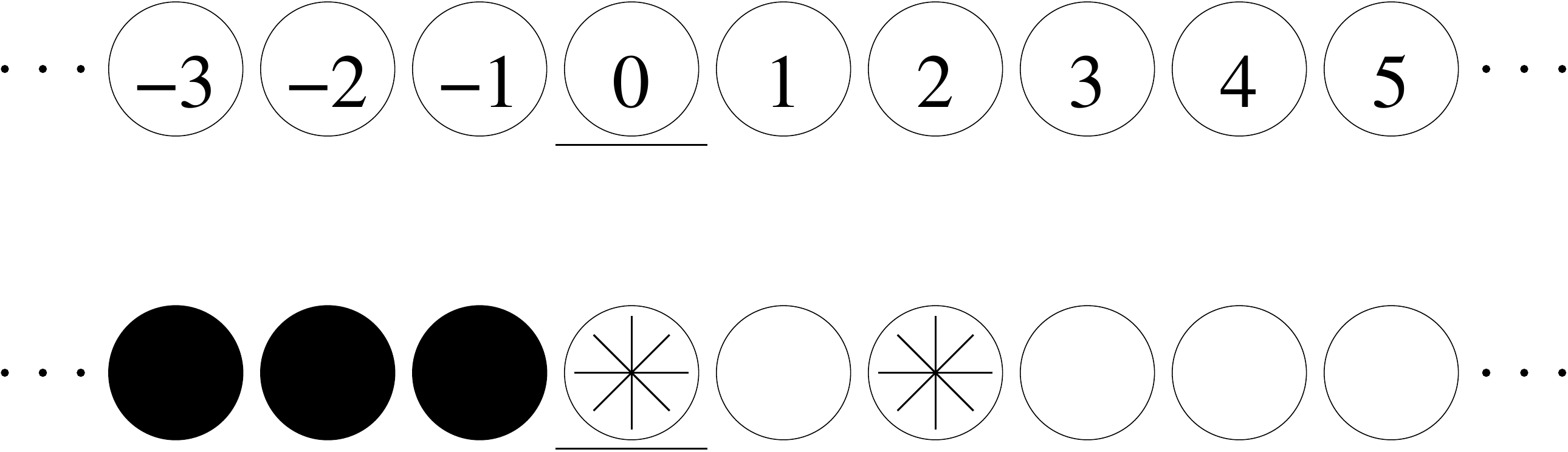}
\caption{\label{initialfig}
Aligned initial conditions for the multi-type \TASEP{} (above) and the two-type \TASEP{} (below)}
\end{figure}

To deduce Corollary \ref{distancetwocorr} from Theorem \ref{t2}(b), one can
compare the multi-type \TASEP{} of this section with the two-type \TASEP{}
of Theorem \ref{t2}(b), started from the initial condition 
${\hat\parti}{\un\second}{\hole}{\second}{\hat\hole}$, as shown in Figure \ref{initialfig}.
\newline
Let the two processes evolve together under the basic coupling; 
potential jumps occur at the same site at the same times in the two processes.

Folllowing the joint evolution of the two processes, one can observe  
that the site occupied by particle 0 (in the multi-type process) always corresponds
to the site of a second-class particle (in the two-type process).
In addition, the other second-class particle is at the site either of particle 1 or of particle 2, whichever 
is further to the right at the time. 
So the first time that one second-class particle tries to jump onto the site
of the other is precisely the time that particle 0 jumps onto the site of whichever is 
further right of particle 1 and particle 2, i.e.\ the first moment that particle 0 has
overtaken both particle 1 and particle 2. 

A similar argument can be used to show the following; 
the probability that particle $i$ overtakes 
all of particles $i+1, i+2,\dots, i+m$ in the multi-type process
is the same as the probability that the two second-class particles
\coalesce in a two-type process starting from the state
${\hat\parti}{\un\second}{\hole\dots\hole}{\second}{\hat\hole}$,
where there is a string of $m-1$ holes between the two second-class particles.
We have seen that this probability is $2/3$ for $m=1$ and $1/2$ for $m=2$.
Simulations suggest that the probability is in fact $2/(m+2)$ for all $m$.

\section{\TASEP{} and competition growth}
\label{growthsec}
In this section we describe the two-type growth model which corresponds to 
the \TASEP{} with two second-class particles.

The connection between the \TASEP{} and the corner growth model dates back 
to Rost \cite{Rost}. In the simplest case, the growth
takes place in the positive integer quadrant $\NN^2$ 
(where we write $\NN=\{1,2,\dots\}$).
At time 0 all the sites of this quadrant are empty.
A site $(x,y)$ may become occupied if both sites $(x-1,y)$ 
and $(x,y-1)$ are occupied; once this condition holds, 
the time until site $(x,y)$ becomes occupied is an
exponential random variable with rate 1 (independently for each site).
As a boundary condition one assumes that the sites $(0,x)$ and $(x,0)$ are 
always occupied, for each $x\in\NN$.

We will consider an extension of this model, where each occupied
site has one of two colours, red and blue. A site 
takes its colour at the moment it is occupied and never changes thereafter. 
In the initial condition, 
the sites $(1,1)$,
$(1,2)$ and $(2,1)$, and all other sites are unoccupied.  (The colour assigned
to the three already-occupied points is not important).

When a site becomes occupied, its colour is chosen as follows:
\begin{itemize}
\item Any site $(1,x)$ or $(x,1)$, for $x>1$, becomes blue.
\item The site $(2,2)$ becomes red.
\item A site $(x,y)$, $x>1$, $y>1$, $(x,y)\ne(2,2)$, takes the colour of whichever
of its parents became occupied \textit{more recently}.
\end{itemize}
(Since the delay times for occupation of each site are independent continuous
random variables, with probability one all times of occupation are different and
so there are no ties to break in the third case above).

There are three possibilities for the first event of the growth process:
$(1,3)$ becomes occupied and blue, or $(3,1)$ becomes occupied and blue,
or $(2,2)$ becomes occupied and red. Each of these occurs with probability $1/3$.

\begin{figure}[h]
\begin{minipage}[t]{0.25\textwidth}
\includegraphics[width=0.99\linewidth]{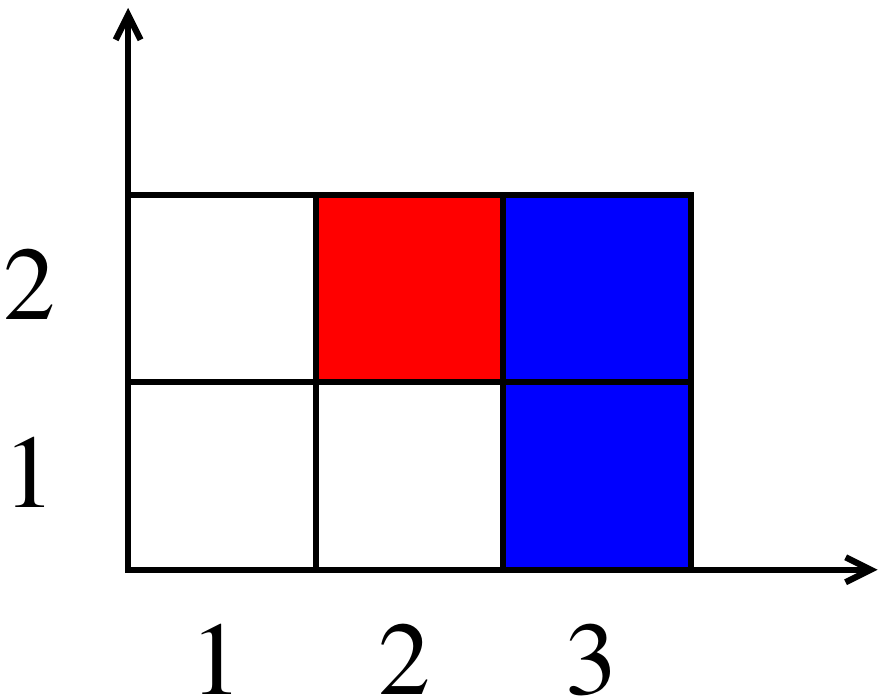}
\end{minipage}
\qquad\qquad
\begin{minipage}[t]{0.25\textwidth}
\includegraphics[width=0.99\linewidth]{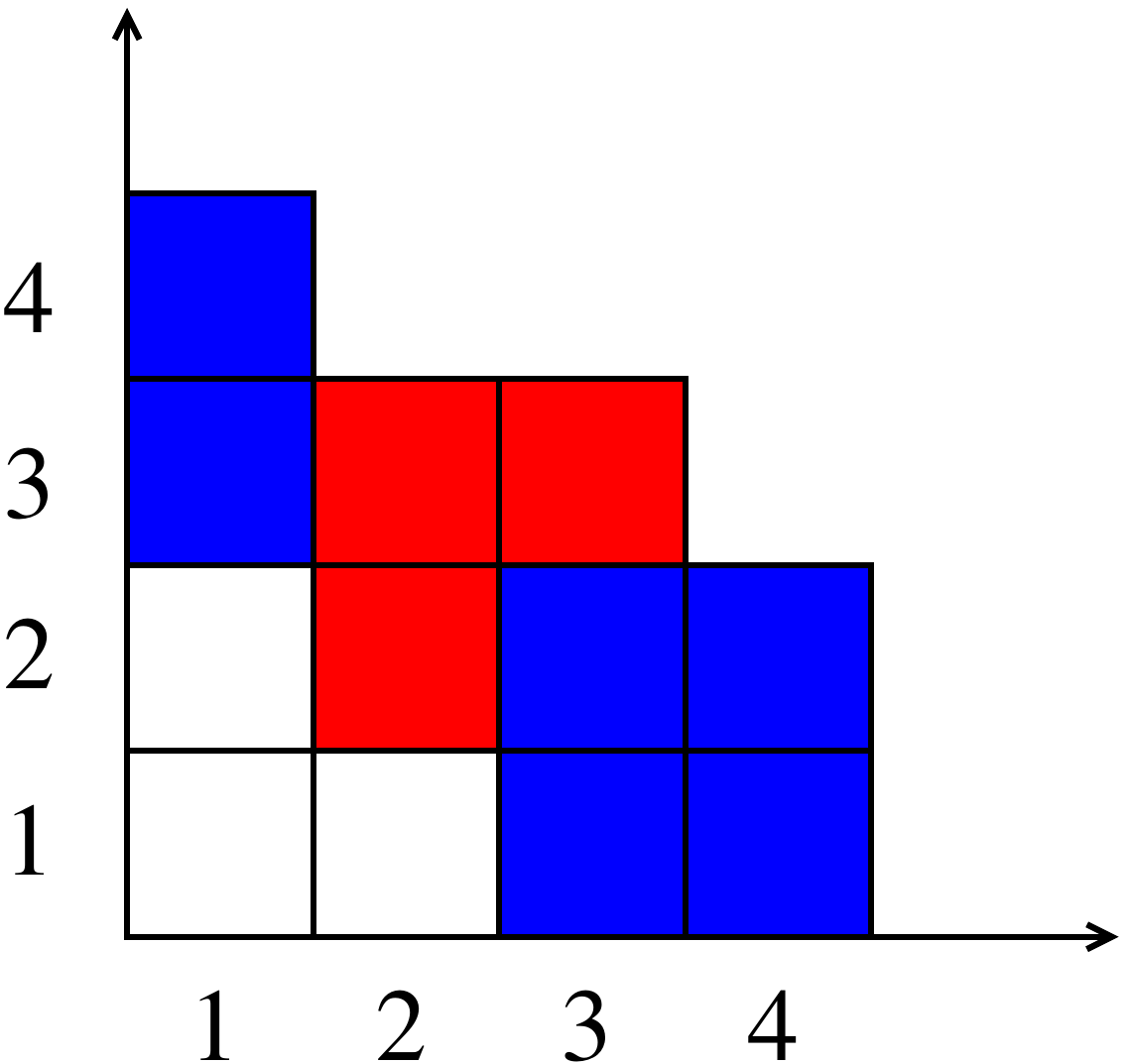}
\end{minipage}
\caption{\label{examplefig} To the left an example showing a growth process 
after its first three jumps.  We can see that $G(2,2)<G(3,1)$, since
  site $(3,2)$ took its colour from $(3,1)$ rather than $(2,2)$. Also clearly
  $G(2,2)<G(1,3)$ (since $(1,3)$ is not yet occupied).  As a result site $(2,3)$
  will take its colour from $(1,3)$ and so will also be blue; the red cluster
  will be surrounded and unable to grow further. To the right, an example
  showing a growth process after its first nine jumps. 
Here one can deduce, for example, that the first site to be 
occupied was $(1,3)$. Given that $(1,3)$ is the first to be occupied, 
the conditional probability that the red cluster survives is 1/2 (see 
the remark at the end of this section).}
\end{figure}

Note that there will certainly be infinitely many points which 
take the colour blue (for example, all the points on the axes $(x,1)$ and $(1,x)$).
However, there may or may not be infinitely many red points.

For example, suppose that the first event in the process is 
that $(2,2)$ is occupied (red). In this case, both $(1,3)$ and
$(3,1)$ are occupied later than $(2,2)$. As a result,
$(2,3)$ will take its colour from $(1,3)$ (its more recently
occupied parent) and similarly $(3,2)$ will take its colour from $(3,1)$.
So both $(3,1)$ and $(1,3)$ become blue, and the red cluster has been
``surrounded''; there will be no further red sites. 

If $(2,2)$ is not the first site to become occupied, 
the red cluster may nevertheless become surrounded later;
alternatively, it may grow unboundedly
See Figures \ref{examplefig} 
and \ref{meetingfig} for some examples.




\begin{theorem}\label{growththm}
The probability that the red cluster grows unboundedly is $1/3$. 
\end{theorem}

\begin{figure}[h]
\begin{minipage}[t]{0.36\textwidth}
\centering
\includegraphics[width=0.99\linewidth]{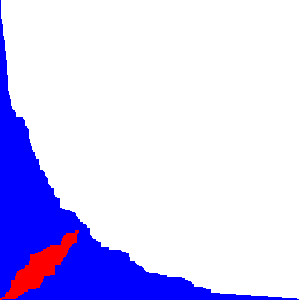}
\end{minipage}
\begin{minipage}[t]{0.36\textwidth}
\centering
\includegraphics[width=0.99\linewidth]{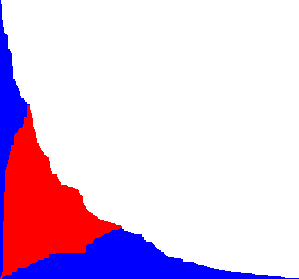}
\end{minipage}
\caption{\label{threeexamplesfig}
  Two simulations of the growth process up to time 300.  The left
  simulation is a rare example where the red cluster grows for a long time (here
  up until time 271.8) but then becomes surrounded. On the right, the red
  cluster looks likely to survive for ever.}
\end{figure}

Coupier and Heinrich \cite{CouHei} recently showed that the probability in
Theorem \ref{growththm} was strictly positive, without using any \TASEP{}
representation, by proving that, in the associated directed last-passage
percolation model, there exists with positive probability an infinite geodesic
starting from the point $(1,1)$ and passing through the point $(2,2)$.

Theorem \ref{growththm} is again deduced from Theorem \ref{t1}(a).  We give the
details of the proof in Section \ref{growthproofsec}.  The main element in the
proof of Theorem \ref{growththm} is the representation introduced by Ferrari and
Pimentel \cite{F.P.} in which a site containing a second-class particle particle
$\second$ is replaced by a pair of sites in the state $[\hole, \parti]$.
Ferrari and Pimentel treated systems with just one second-class particle, so
that the two-type system with initial condition ${\hat\parti}\second{\hat\hole}$
was replaced by a one-type system with initial condition
${\hat\parti}\hole\parti{\hat\hole}$.  In Rost's correspondence, this equates to
a growth process starting from the state where site $(1,1)$ alone is occupied.

Analogously, the initial condition ${\hat\parti}\second\second{\hat\hole}$ from
Theorem \ref{t1}(a) corresponds to the initial condition
${\hat\parti}\hole\parti\hole\parti{\hat\hole}$, 
or to a growth model
starting with sites $(1,1)$, $(2,1)$ and $(1,2)$ occupied. 

This hole-particle representation of two second-classs particles 
survives only until the moment when one second-class
particle attempts to jump onto another; 
but this is enough for our purposes, since we will be 
concerned only with the probability of such a collision event occurring.

In \cite{F.P.}, the path of one second-class particle corresponded to the 
interface between two clusters. Here we have two interfaces, between 
the red cluster in the middle and the blue clusters on either side of it.
The event that the two second-class particles meet will correspond to the event 
the two interfaces meet. If this happens, the red cluster between the two interfaces is
surrounded and is unable to grow any further. This will happen 
with probability 2/3, as in Theorem \ref{t1}(a). 
See Section \ref{growthproofsec} for full details.

\paragraph{\sl Remark}
Note that the first site to be occupied in the growth process is 
equally likely to be $(1,3)$, $(2,2)$ or $(3,1)$, with probability
$1/3$ each. As observed above, if site $(2,2)$ is the first to be occupied,
the red cluster cannot survive. 
We know that the overall probability that the red cluster survives is $1/3$.
Hence, 
given that the first site to be occupied is \textit{not} $(2,2)$, 
the conditional probability that the red cluster survives is $1/2$.
This statement corresponds to part (b) of Theorem \ref{t1}, in the case
$p=1$, (while 
Theorem \ref{growththm} corresponds to part (a) of Theorem \ref{t1}).
By symmetry, it does not matter whether $(3,1)$ or $(1,3)$ is 
occupied first.
To emphasise the connection, note that if we start from the configuration
where $(1,3)$ has been occupied in addition to $(1,1)$, $(1,2)$ and $(2,1)$,
this corresponds to starting the \TASEP from the configuration 
$\hatparti\hole\parti\hole\hole\parti\hathole$, 
which, under the hole-particle representation of the second-class particle,
corresponds in turn to the initial condition 
$\hatparti\second\hole\second\hathole$ of Theorem \ref{t1}(b).
Again, see Section \ref{growthproofsec} for details 
of the correspondence between states of the \TASEP and states of the 
growth model. 

\section{Hydrodynamics and coupling}
\label{s3}

In this section we state some known results on hydrodynamics
and on coupling that will be used in the proofs.

\paragraph{\bf Hydrodynamics}
Define the product measures
$\nu_\lambda, 0\leq\lambda\leq 1$
under which every site is occupied by a particle with probability
$\lambda$ and a hole with probability $1-\lambda$, independently for different sites.
Any translation-invariant stationary distribution of the one-type \ASEP{}
is a mixture of these product measures
\cite{L.}.

Suppose the initial distribution puts probability 1 on configurations with 
asymptotic density 1 to the left and 0 to the right. (In particular
this is the case whenever all sites to the left of some point contain 
particles and all sites to the right of some point contain holes, i.e.\
the initial configuration has the form $\hat\parti\, a\, \hat\hole$ for some finite string $a$).

Then the process converges to local equilibrium:
\begin{equation}\label{localequilibrium}
  \lim_{N\rightarrow{\infty}}S(tN)\theta_{\lfloor rN\rfloor}f=\nu_{u(r,t)}f
\end{equation}
for any cylinder function $f$, where $\theta_{x}$ is the translation by $x$,
$S(t)$ is the semigroup of the \ASEP{} with generator given by (\ref{a30gen}),
$\lfloor\cdot\rfloor$ is the integer part and for $t\geq{0}$, $r\in{\mathbb{R}}$,
$u(r,t)$ is the entropy solution of the inviscid Burgers equation
\begin{equation*}
\left\{
  \begin{array}{ll}
    \partial_{t} u(r,t)+(p-q)\partial_{r} \Big(u(r,t)(1-u(r,t))\Big)=0\\
    u(r,0)=I(r<0)
  \end{array},
\right.
\end{equation*}
which is explicitly given by
\begin{equation}
\label{a38}
u(r,t)\,=\,
\left\{
         \begin{array}{ll}
           1, & \hbox{$r\leq{-(p-q)t}$;} \\
           \frac{(p-q) t-r}{2(p-q) t},
           & \hbox{$-(p-q)t<r\leq{(p-q)t}$;} \\
           0, & \hbox{$r>(p-q)t$.}
         \end{array}
       \right.
\end{equation}

\paragraph{\bf Basic Coupling} The Harris construction of the process makes use
of a family of independent Poisson processes $\omega=(N(x,x\pm1),\,x\in\mathbb
Z)$, where $N(x,x+1)$ has rate $p$ and $N(x,x-1)$ has rate $q$ .
These live in a probability space
$(\Omega,\mathcal F, P)$. 
For the process with first-class and second-class particles, 
we write $\Xi=\{\hole,\second,\parti\}^{\mathbb Z}$ as above, 
and we define a function $\Phi:\Omega\times\Xi\times
R^+\to \Xi$ inductively in time as follows (informally): fix an initial
configuration $\eta$ and call $\eta_t=\Phi(\omega, \eta,t)$. Then
\[
\eta_t:= \left\{\begin{array}{ll}
    (\eta_{t-})^{x,y}, & \hbox{if } t\in N(x,y) \hbox{ and }\eta(x)\prec\eta(y) \\
\eta_{t-},& \hbox{otherwise}
  \end{array}\right..
\]
that is, the particles attempt to jump from $x$ to $y$ at the epochs of the
Poisson process $N(x,y)$ and do so if their classes admit. Of course this
definition is incomplete, as the epochs of the Poisson process are not well
ordered, but a simple argument shows that for all $t$ and almost all realizations
$\omega$, $\mathbb Z$ is partitioned into finite intervals that do not interact with
each other in the interval $[0,t]$ (see \cite{FM1} and \cite{FM2} for details). The
process so defined is Markov and has generator $L$ given in \eqref{a30gen}.

The basic coupling between $n$ versions of the process with initial
configurations $\un{\eta}=(\eta^1,\dots,\eta^n)\in \Xi^n$ is the process which
uses the same $\omega$ for all marginals:
\begin{equation}
  \label{bc1}
  \un{\eta}_t=(\eta^1_t,\dots,\eta^n_t)
  := \big(\Phi(\omega,\eta^1,t),\dots,\Phi(\omega,\eta^n,t)\big)
\end{equation}
That is, in the basic coupling the particles attemp to jump from $x$ to $y$ at
the same epochs in each marginal, always respecting the classes (at the
marginal).

\section{Speed of a second-class particle}
\label{ss4}

\noindent\textbf{Proof of Theorem \ref{2cptheorem}}.
 As in the proof of Theorem 1 of \cite{F.K.}, the main ingredient is a coupling
 argument and the convergence to local equilibrium (\ref{localequilibrium}).
 One wants to show that for $r\in[-(p-q), p-q]$,
\begin{equation}\label{form}
  \lim_{t\rightarrow{\infty}}P\big(X(t)>rt\big)
=\frac{p-q-r}{2(p-q)}.
\end{equation}
We will consider \ASEP{} processes $\eta_t$ and $\zeta_t$
started from the initial states $\eta_0={\hat\parti}\un\parti{\hat\hole}$.
and $\zeta_0={\hat\parti}\un\hole{\hat\hole}$.

Let $J_t^r(\eta_0)$ be the number of particles of $\eta_t$ which
are to the right of $rt$, and similarly for $J_t^r(\zeta_0)$.

We will compute $E (J_t^r(\eta_0))-E (J_t^r(\zeta_0))$ in two different ways.

The initial states $\eta_0$ and $\zeta_0$ differ only at the origin,
where $\eta_0$ has an extra particle.
If we couple the processes using the basic coupling, then there
is exactly one discrepancy at all times, and it
behaves as the single second-class particle of a
process started in the state ${\hat\parti}\un\second{\hat\hole}$.
Hence
\begin{equation}
\label{way1} E (J_t^r(\eta_0))-E (J_t^r(\zeta_0))=P\big(X(t)>rt\big).
\end{equation}

On the other hand, by shifting the origin by 1,
we can see that the number of particles to the right of $rt$
in $\zeta_t$ has the same distribution as the number of particles
to the right of $rt+1$ in $\eta_t$. Hence
\begin{equation}\label{way2}
E (J_t^r(\eta_0))-E (J_t^r(\zeta_0))= P\Big(\eta_t\left(\lceil rt \rceil\right)=1\Big),
\end{equation}
which, by the convergence to local equilibrium in
(\ref{localequilibrium}), is $\frac{p-q-r}{2(p-q)}$,
so that (\ref{way1}) and (\ref{way2}) together
give (\ref{form}) as desired.
\cqd

\section{Distance between two second-class particles}
\label{ss5}

In this section we give  the proof of Theorem \ref{t2}. 
For $r\in[-1,1]$ let
\begin{equation}
  \label{c1}
  J_{t}^{r}(\eta)\,:=\,\sum_{x>rt}\eta_{t}(x) + \eta_t(\lfloor rt\rfloor)  (\lfloor rt\rfloor+1-rt)
\end{equation}
 be the number of particles at the right of $rt$ at time $t$; for
convenience, we consider that if there is a particle at $\lfloor rt\rfloor$ a fraction
$(\lfloor rt\rfloor+1-rt)$ of it is to the right of $rt$. 
\medskip

\paragraph{\bf Proof of Theorem \ref{t2}}

Consider an \ASEP{} $\eta_t, t\geq0$
started from the initial condition $\eta_0=
\hat\eta=
{\hat\parti}\un\parti{\hat\hole}
=
{\hat\parti}\un\parti\hole{\hat\hole}$.

From $\hat\eta$, the first jump of the process is to the state 
${\hat\hole}\un\hole\parti{\hat\hole}$
and this jump occurs
at rate $p$.

Hence using the
definition of $J_t^r$ in (\ref{c1}) and
the Kolmogorov backward equation,
we have that $EJ_t^r(\hat\eta)$ is continuous
and is differentiable everywhere except possibly at points $t$ where $rt\in\mathbb{Z}$, with
derivative
\begin{equation}
  \label{current2}
  \frac{d}{dt}E(J_{t}^{r}(\hat\eta)
=p\Big[E\Big(J_{t}^{r} \big({\hat\parti}\un\hole\parti{\hat\hole} \big)
-E(J_{t}^{r}\big(  {\hat\parti}\un\parti\hole{\hat\hole}      \big)\Big)\Big]
-rE(\eta_t(\lfloor rt\rfloor)).
\end{equation}
We couple two
versions of the \ASEP{} starting from ${\hat\parti}\un\hole\parti{\hat\hole}$ and from
${\hat\parti}\un\parti\hole{\hat\hole}$, aligned as follows:
\begin{equation}
\begin{matrix}
{\hat\parti}\un\hole\parti{\hat\hole} \\ 
{\hat\parti}\un\parti\hole{\hat\hole} 
\end{matrix}
\end{equation}
At time zero there are two discrepancies, located at sites 0 and 1. If a jump is 
attempted from the site of one the discrepancies into the
other, the discrepancies disappear and the two configurations coincide from that 
instant on. Until that time, denoted by $\tau$, the discrepancies behave as the positions of 
two second-class particles, $X(t)$ and $Y(t)$.
The difference of the number of particles to the right of $rt$ for the two coupled 
processes can be expressed in terms of the positions of $X(t)$
and $Y(t)$, as follows:
\begin{equation}
\label{c2}
J^r_{t}\big(  {\hat\parti}\un\hole\parti{\hat\hole}  )-J^r_{t}(   {\hat\parti}\un\parti\hole{\hat\hole}    \big)
=
\begin{cases}
1, & \hbox{if $\tau>t$ and $X(t)<\lfloor rt\rfloor<Y(t)$;} \\
rt-\lfloor rt\rfloor, & \hbox{if $\tau>t$ and $X(t)=\lfloor rt\rfloor<Y(t)$;} \\
1-rt+\lfloor rt\rfloor, & \hbox{if $\tau>t$ and $X(t)<\lfloor rt\rfloor=Y(t)$;} \\
0,& \hbox{otherwise.}
\end{cases}
\end{equation}
Inserting \eqref{c2} in (\ref{current2}),
\begin{multline}
  \label{a37}
  \frac{d}{dt}E(J^r_{t}(\hat\eta))
  =
p\Big[
P\big(X(t)<\lfloor rt\rfloor<Y(t)\big)
+(rt-\lfloor rt\rfloor) P\big(X(t)=\lfloor rt\rfloor<Y(t)\big)
\\
+(1-rt+\lfloor rt\rfloor) P\big(X(t)<\lfloor rt\rfloor=Y(t)\big)
\Big]
- rE(\eta_t(\lfloor rt\rfloor))
\end{multline}

Now we compute (\ref{current2}) in a different way by associating to
$J_{t}^{r}(\eta_0)$ the following martingale:
\begin{equation*}
  M^r_{t}
  =J_{t}^{r}(\eta_0)
-\int_{0}^{t}
\left\{
(\lfloor rs\rfloor+1-rs)w_{\lfloor rs\rfloor-1}(\eta_s)
+
(rs-\lfloor rs\rfloor)w_{\lfloor rs\rfloor}(\eta_s)
-r \eta_{s}(\lfloor rs\rfloor)
\right\}\, ds,
\end{equation*}
where, for a fixed configuration $\eta$, $w_{x}(\eta)$ denotes the
instantaneous current through the bond $[x,x+1]$, namely:
\begin{equation*}
  w_{x}(\eta)\,=\,p\,\eta(x)(1-\eta(x+1))\,-\,q\,\eta(x+1)(1-\eta(x)).
\end{equation*}
This expression together with the convergence to local equilibrium
(\ref{localequilibrium}) imply
\begin{equation}
\label{a41} \lim_{t\rightarrow{\infty}}
\frac{d}{dt}E(J^r_{t}(\hat\eta))=(p-q)u(r,1)(1-u(r,1))-ru(r,1).
\end{equation}
The same convergence to local equilibrium also gives that
\begin{equation}\label{easyconv}
E(\eta_t(\lfloor rt\rfloor))\to u(r,1) \text{ as } t\to\infty.
\end{equation}
Note that
\begin{equation}\label{fracsimp}
\lim_{t\to\infty}P\big(X(t)=\lfloor rt\rfloor\big)=\lim_{t\to\infty}P\big(Y(t)=\lfloor rt\rfloor\big)=0.
\end{equation}
(For example, $X(t)$ represents the discrepancy between 
two coupled processes starting from the states
${\hat\parti}\un\hole{\hat\hole}$ 
and ${\hat\parti}\un\parti{\hat\hole}$; the probability of finding this discrepancy
precisely at the site $\lfloor rt\rfloor$ at time $t$ goes to 0 as $t\to\infty$ since,
again by the convergence to local equilibrium in (\ref{localequilibrium}),
the probability of site $\lfloor rt\rfloor$ being occupied at time $t$ 
converges to the same limit $u(r,1)$ for both processes.)

Putting equations (\ref{a37})-(\ref{fracsimp}) together gives
\[
\lim_{t\to\infty} P\big(X(t)\leq rt<Y(t)\big)=\frac{p-q}{p}u(r,1)\big[1-(u(r,1))\big].
\]
The form of the function $u$ from \eqref{a38} then gives (\ref{separationprob}).

The fact that $X(t)<Y(t)$ whenever $\tau>t$, and the relation
\[
E\Big[\big(Y(t)-X(t)\big)^+I(\tau>t)\Big]=
\sum_y P\Big(X(t)\leq y<Y(t), \tau>t\Big),
\]
then give
\[
\lim_{t\to\infty}\frac{1}{t}E\Big[\big(Y(t)-X(t)\big)I(\tau>t)\Big]=
\frac{1}{4p(p-q)}\int_{-(p-q)}^{p-q}\left[(p-q)^2-r^2\right]dr,
\]
giving \eqref{eq1} as required.\cqd

\section{Collision probabilities}
\label{ss6}

\paragraph{\bf Proof of Theorem \ref{t1}(a)}

Denote by 
$\hatparti\un\second\hole\hathole$
the configuration that has a second-class
particle at the origin, while the negative sites are occupied by first-class
particles and the rest is empty. Denote by 
$\hatparti\un\parti\second\hathole$
the configuration that has a second-class particle at site 1, while the non-positive
sites are occupied by first-class particles and the rest is empty. We couple
both processes and prove later that the positions of the discrepancies initially
at sites $0$ and $1$ behave as $X(t)$ and $Y(t)$, the two second-class
particles of the theorem, with $X(0)=0$ and $Y(0)=1$.

For a configuration $\eta$ with a unique second-class particle and a finite
number of first-class particles to the right of it, denote by $X(t,\eta)$ the
position of the second-class particle at time $t$ and by $J^{(2)}_{t}(\eta)$ the
number of first-class particles to the right of $X(t,\eta)$ at time $t$:
\begin{equation*}
  J_{t}^{(2)}(\eta):=\sum_{x\geq{1}}\eta_{t}\big(X(t,\eta)+x\big).
\end{equation*}
where $\eta_t$ is the configuration at time $t$ for the process with initial
configuration $\eta_0=\eta$.

Define the configuration
$\teta=
\hatparti\un\second\hathole$. For $\eta_0=\teta$,
the Kolmogorov backwards equation gives
\begin{equation}\label{a22}
  \frac{d}{dt}E\Big(J_{t}^{(2)}(\teta)\Big)
  =
p E\Big(J_{t}^{(2)}(\hatparti\second\un\parti\hathole )\Big)
+p E\Big(J_{t}^{(2)}(\hatparti\un\hole\second\hathole)\Big)
-2pE\Big(J_{t}^{(2)}(\teta)\Big)\Big).
\end{equation}
In \eqref{a22} $\hatparti\second\un\parti\hathole$      
is attained when the rightmost first-class
particle at site $-1$ jumps to $0$, the site occupied by the second-class
particle,
and $\hatparti\un\hole\second\hathole$
arises when the second-class particle jumps
from $0$ to $1$, site occupied by the leftmost hole. The law of
$J_{t}^{(2)}(\eta)$ does not depend on the actual location of $X(0,\eta)$ but on
the relative positions of the other particles with respect to $X(0,\eta)$. So,
we are free to change the origin. In particular
$E\Big(J_{t}^{(2)}(\hatparti\second\un{\parti}\hathole)\Big)=
E\Big(J_{t}^{(2)}(\hatparti\un{\second}\parti\hathole)\Big)$ and
$2E\Big(J_{t}^{(2)}(\hatparti\un{\second}\hathole)\Big)=
E\Big(J_{t}^{(2)}(\hatparti\un{\second}\hole\hathole)\Big)+E\Big(J_{t}^{(2)}(\hatparti\un{\parti}\second\hathole)\Big)$.
Using this, for any (four-marginals) coupling, we can rewrite \eqref{a22} as
\begin{equation}
\label{a21}
\frac{d}{dt}E\Big(J_{t}^{(2)}(\teta)\Big)
=pE
\big[
J_{t}^{(2)}( \hatparti\un{\second}\parti\hathole )
-J_{t}^{(2)}( \hatparti\un{\second}\hole\hathole  )
+J_{t}^{(2)}( \hatparti\un\hole\second\hathole     )
-J_{t}^{(2)}( \hatparti\un{\parti}\second\hathole   )
\big],
\end{equation}
so that we are dealing with four initial configurations that are aligned as follows:
\begin{equation}
  \label{a23}
\begin{matrix}
 \hatparti\un{\second}\parti\hathole \\
\hatparti\un{\second}\hole\hathole\\  
\hatparti\un\hole\second\hathole    \\ 
\hatparti\un{\parti}\second\hathole   
\end{matrix}
\end{equation}
We perform the basic coupling, where recall the jumps are attempted at the same
time at the same sites for the four marginals. The first and third lines
contribute with a plus sign to \eqref{a21} while the others do so with a minus
sign. We start with two discrepancies at the sites 0 and 1. As before, under the basic coupling
each of these discrepancies moves as a second-class particle, 
up until the moment when a jump is attempted from the site of one onto the site of the other. 
Let $\tau$ denote this time, which could be infinity. 
The states of the four marginals at various times look as follows:
\newcommand{\qq}{\qquad}
\begin{equation}
  \label{a25}
  \begin{matrix}
\hbox{Time 0}&\qq&\hbox{Time $t\in(0,\tau)$}&\qq&\hbox{Time $\tau-$}&&&\hbox{\hspace{-1.3cm}Time $\tau$\hspace{-1.3cm}}&\\
\hatparti \un{\second }\parti \hathole &\qq&\hatparti a\second b\parti c\hathole &\qq&\hatparti a\second \parti c\hathole &\qq&\hatparti a\second \parti c\hathole &&\hatparti a\parti \second c\hathole\\
\hatparti \un{\second }\hole \hathole &\qq&\hatparti a\second b\hole c\hathole &\qq&\hatparti a\second \hole c\hathole &\qq&\hatparti a\hole \second c\hathole &\textbf{or}&\hatparti a \second \hole c\hathole\\
\hatparti \un{\hole }\second \hathole &\qq&\hatparti a\hole b\second c\hathole &\qq&\hatparti a\hole \second c\hathole &\qq&\hatparti a\hole \second c\hathole &&\hatparti a\second \hole c\hathole\\
\hatparti \un{\parti }\second \hathole &\qq&\hatparti a\parti b\second c\hathole &\qq&\hatparti a\parti \second c\hathole &\qq&\hatparti a\second \parti c\hathole &&\hatparti a\parti \second c\hathole
\end{matrix}
\end{equation}
Here $a$, $b$ and $c$ are some finite sequences of holes and first-class particles.
(They are the same within any given column, though
their values may of course differ between columns).
The two possibilities corresponding to time $\tau$ reflect the fact that 
the jump at time $\tau$ may be from the discrepancy on the left onto the discrepancy 
on the right, or vice versa. If the jump is from left to right, one gets the left column; 
if the jump is from right to left, one gets the right column. 

Note that before time $\tau$, the first row has one more $\parti$ to the
right of the $\second$ than the second row does.
Meanwhile the third and the fourth rows have the same number of particles
to the right of the $\second$.
Hence for any $t<\tau$,
\[
J_{t}^{(2)}(\hatparti\un{\second}\parti\hathole)
-J_{t}^{(2)}(\hatparti\un{\second}\hole\hathole)
+J_{t}^{(2)}(\hatparti\un\hole\second\hathole)
-J_{t}^{(2)}(\hatparti\un{\parti}\second\hathole)
=1.
\]

At time $\tau$, the first and the fourth rows become equal,
and (since they are coupled) will still be equal at all later times.
The same is true of the second and third rows.
Hence for any $t\geq\tau$,
\[
J_{t}^{(2)}(\hatparti\un{\second}\parti\hathole)
-J_{t}^{(2)}(\hatparti\un{\second}\hole\hathole)
+J_{t}^{(2)}(\hatparti\un\hole\second\hathole)
-J_{t}^{(2)}(\hatparti\un{\parti}\second\hathole)
=0.
\]
Thus we have, for any $t$,
\[
J_{t}^{(2)}(\hatparti\un{\second}\parti\hathole)
-J_{t}^{(2)}(\hatparti\un{\second}\hole\hathole)
+J_{t}^{(2)}(\hatparti\un\hole\second\hathole)
-J_{t}^{(2)}(\hatparti\un{\parti}\second\hathole)
=I(t<\tau).
\]

Substituting into \eqref{a21}, we get
\begin{equation}
  \label{a32}
  \frac{d}{dt}E(J_{t}^{(2)}(\teta))=
p P(t<\tau)
\end{equation}
Since $\{s<\tau\}$ are decreasing events in $s$, the following limits exist:
\begin{equation}
\lim_{t\rightarrow{\infty}}P(t<\tau)
\,=\,\lim_{t\rightarrow{\infty}}\frac{1}{t}\int_{0}^{t}P(s<\tau)ds
\,=\, \lim_{t\rightarrow{\infty}}\frac{1}{t}\int_{0}^{t}\frac{1}{p}\frac{d}{ds}E(J_{s}^{(2)}(\teta))ds
\,=\, \frac{1}{p}\lim_{t\rightarrow{\infty}}\frac{E(J_{t}^{(2)}(\teta))}{t}.
\end{equation}

By Theorem \ref{2cptheorem}, $X(\teta)/t$ converges in 
distribution as $t\rightarrow{\infty}$ to a Uniform random variable with
support on $[-(p-q),(p-q)]$ denoted again by $\mathcal{U}_{p}$. 
Together with the convergence to equilibrium in (\ref{a38}), this allows us to
conclude that
\begin{equation}
\label{convergence of current in distribution}
  \frac{J_{t}^{(2)}(\teta)}{t}\xrightarrow[t\rightarrow{\infty}]
  \,\frac{1}{p-q}\Big(\frac{(p-q)-\mathcal{U}}{2}\Big)^2
\end{equation}
in distribution.
In addition, $J_{t}^{(2)}(\teta)$ is non-negative and dominated
by a Poisson process of rate 1. Hence
\begin{equation}\label{a48}
  \lim_{t\rightarrow{\infty}}E\Big(\frac{J^{(2)}_{t}(\teta)}{t}\Big)
  =\frac{1}{p-q}E\Big[\Big(\frac{(p-q)-\mathcal{U}}{2}\Big)^2\Big]
= \frac{p-q}3.
\end{equation}
Taking the limit as $t\to\infty$ in (\ref{a32}) gives
$P(\tau<\infty)=1-\frac{p-q}{3p}=\frac{1+p}{3p}$ as required.
\cqd

\noindent\textbf{Proof of Corollary \ref{2cpcorollary}}.
Using \eqref{a21} and the forward equations, we have
\begin{align}
\nonumber
p P(t<\tau)
&=\frac{d}{dt}E(J_{t}^{(2)}(\teta))
\\
\label{corruse1}
&=\,pE\eta_{t}\big(X(t)-1)\big)
-qE\eta_{t}\big(X(t)+1\big),
\end{align}
since $J_{t}^{(2)}$ increases by one whenever the second-class particle
swaps with a particle on its left, and decreases by one whenever
the second-class particle swaps with a particle on its right.

On the other hand, by the symmetry of the system under
exchange of hole and particle and of left and right,
$\eta_{t}\big(X(t)-1\big)$ has the same law as $1-\eta_{t}\big(X(t)+1\big)$. 
Hence,
\begin{equation}
\label{corruse2}
E\eta_{t}\big(X(t)-1\big)
+E\eta_{t}\big(X(t)+1)\big) = 1.
\end{equation}

As $t\to\infty$, the left-hand side of \eqref{corruse1}
converges to $pP(\tau=\infty)$ which is $(p-q)/3$ by Theorem \ref{t1}.
Putting \eqref{corruse1} and \eqref{corruse2} together then gives \eqref{correq1} and 
\eqref{correq2} as desired.
\cqd

\noindent\textbf{Proof of Theorem \ref{t1}(b)}.
Consider the process started from the state $\hatparti\un\second\second\hathole$.
From Theorem \ref{t1}(a), we know that the probability of collision
of the two second-class particles is $(1+p)/3p$. 

We condition on the first jump. There are three possibilities. At rate 1, the 
two second-class particles collide (either because of a jump right from site 0 
or because of a jump left from site 1). At rate $p$ there is a jump right by the 
first-class particle at site $-1$, displacing the second-class particle at site 0 and 
leading to the state $\hatparti\second\un\parti\second\hathole$. 
Also at rate $p$ there is a jump right by the second-class particle at site 1, 
leading to the state $\hatparti\un\second\hole\second\hathole$.

So with probability $1/(1+2p)$ the first jump leads to immediate collision,
and with probability $2p/(1+2p)$ the first jump is to one of the states 
$\hatparti\second\un\parti\second\hathole$ or $\hatparti\un\second\hole\second\hathole$.

Notice that the dynamics of the system are invariant under 
the operation which reverses left and right and also exchanges 
the roles of first-class particles and holes. From this one can see that
the probability of collision of the two second-class particles is the same 
from the two states $\hatparti\second\un\parti\second\hathole$
and $\hatparti\un\second\hole\second\hathole$. Let this probability be $\alpha$.
Conditioning on the first jump, we obtain
\[
\frac{1+p}{3p}=\frac{1}{1+2p}+\frac{2p}{1+2p}\alpha.
\]
Solving this we obtain $\alpha=(1+2p^2)/(6p^2)$ as required.

\section{Growth model proof}\label{growthproofsec}
In this section we explain the correspondence between the 
growth model with three types and the \TASEP{} with two second-class particles,
which leads to Theorem \ref{growththm}.

This is a natural extension of the correspondence between a two-type growth model
and the \TASEP with \textit{one} second-class particle, developed by Ferrari and 
Pimentel in \cite{F.P.}. We start by giving an account of that correspondence; 
for a more detailed description, the reader may like to refer to \cite{F.P.} itself.

First we make explicit the correspondence between the (one-type) corner growth model
and the one-type \TASEP{}, as described by Rost in \cite{Rost}.

In the corner growth model, all sites of $\NN^2$ start off empty at time 0;
each site $z$ waits until its two neighbours below and to the left,
$z-(0,1)$ and $z-(1,0)$, have been occupied; it then
waits for a random further amount of time which is exponential with mean 1,
and then becomes occupied itself. All sites $(0,x)$ and $(x,0)$ are 
already occupied at time 0.

Write $G(z)$ for the time that site $z$ becomes occupied in the corner growth model.
We can write a system of recurrences $G(z)=w(z)+\max\{G\big(z-(0,1)\big), G\big(z-(1,0)\big)\}$,
for $z\in\NN^2$, where we set the boundary conditions $G(i,j)=0$ whenever $i=0$ or $j=0$.
Here the $w(z),z\in\NN^2$ are i.i.d.\ exponential random variables with mean 1.

\begin{figure}[h]
\centering
\includegraphics[width=0.48\linewidth]{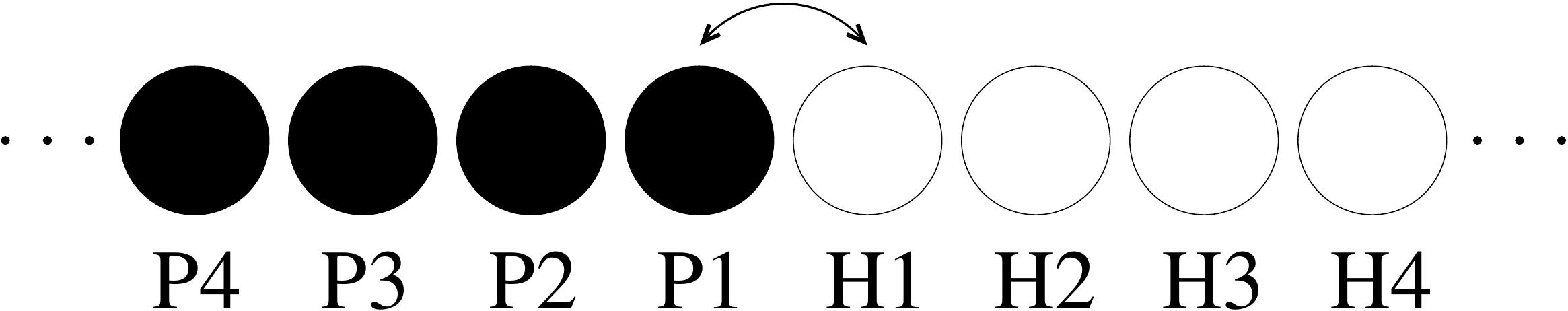}
\caption{\label{sequencefig}
Initial condition for the \TASEP{}, with particles and holes labelled sequentially.
The first event, at time $G(1,1)$, is the exchange of P1 and H1.}
\end{figure}

To represent the \TASEP{} in this way, suppose that the particles are labelled P1, P2, P3, $\ldots$ 
from right to left, and the holes are labelled H1, H2, H3 from left to right,
as shown in Figure \ref{sequencefig}. 
(These labels are not to be confused with priorities or with the types of particle
in Section \ref{growthsec}; 
the sequence of particles and the sequence of holes will both stay in 
their initial orderings throughout, although particles and holes will exchange places).
For $z=(i,j)\in\NN^2$, let $G(i,j)$ be the time at which particle P$i$ overtakes
hole H$j$. Then exactly the same recurrences hold as in the previous paragraph. 
Once particle P$i$ has overtaken hole H$(j-1)$, and hole H$j$ has been overtaken 
by particle P$(i-1)$, particle P$i$ finds itself immediately on the left of hole H$j$;
now an exponential amount of time with rate 1 passes before P$i$ overtakes H$j$.

This gives a one-to-one correspondence between a set of states of the \TASEP{} (in which
all sites far enough left contain particles and all sites far enough right contain holes) and
a set of states of the growth process (with finitely many sites of $\NN^2$ occupied and satisfying the 
condition that if $z$ is occupied then so are $z-(0,1)$ and $z-(1,0)$).
$G(i,j)\leq t$ means that particle P$i$ is to the right of hole H$j$ at time $t$,
or correspondingly that site $(i,j)$ is already occupied at time $t$.

Ferrari and Pimentel \cite{F.P.} extended this correspondence
to represent a \TASEP{} with a single second-class particle in terms of 
two-type growth model.

The first element in this representation was to replace the 
site containing the second-class particle, $\second$, by 
a \textit{pair} of sites, containing a hole and a particle in that order, $[{\hole}{\parti}]$.
This pair is called the \starpair, consisting of the \starhole and the \starparticle.
When it is immediately to the left of a hole, the \starparticle 
jumps onto the hole to its right at rate 1; when this happens, 
we consider that the \starpair 
itself moves one unit to the right, and $[\hole\parti]\hole$ 
becomes $\hole[\hole\parti]$.
Similarly when a particle on the left of the \starpair jumps to the right, 
the \starpair itself moves left, so that $\parti[\hole\parti]$ becomes
$[\hole\parti]\parti$. Hence the \starpair, considered as a single unit, 
behaves precisely as a second-class particle. 

Hence we may represent the initial condition $\hat\parti\second\hat\hole$, 
of a system with a single second-class particle, by the 
initial condition $\hat\parti[\hole\parti]\hat\hole$ of a one-type system. 
Again we label holes and particles sequentially,
and we let $G(i,j)$ be the first time that particle $i$ is to the right of 
hole $j$. The boundary conditions are as before, except that now $G(1,1)=0$ 
(since we start with particle P1 already to the right of hole H1). 

\begin{figure}[h]
\centering
\includegraphics[width=0.6\linewidth]{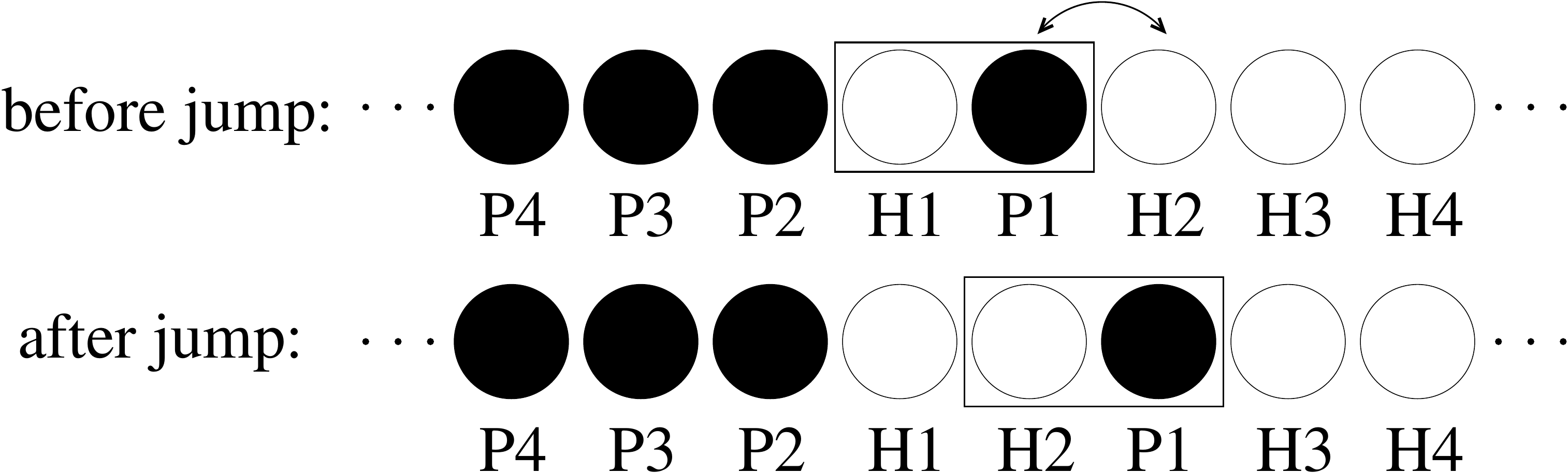}
\caption{\label{jumpexamplefig}Initial condition corresponding
to the system with one second-class particle, and one of the two possible first jumps.
The jump shown occurs at time $G(1,2)$, if $G(1,2)<G(2,1)$.}
\end{figure}

We keep track of the labels of the \starhole and \starparticle at time $t$.
Let these labels be $I(t)$ and $J(t)$. We have $I(0)=J(0)=1$,
since initially the \starpair covers hole H1 and particle P1.
Thereafter, whenever the \starpair jumps right (as in Figure \ref{jumpexamplefig}
for example), $I(t)$ increases by one while $J(t)$ remains the same;
when the \starpair jumps left, $I(t)$ remains the same and $J(t)$
increases by one. In particular, the position $X(t)$ 
of the second-class particle at time $t$ corresponds to $I(t)-J(t)$
(the number of its jumps to the right minus the number of its jumps to the left).
See Figure \ref{jumpexamplefig} for an example. 

The evolution of the process $(I(t),J(t))$ can be described as follows. 
At time $t$, the hole with label $I(t)$ is immediately to the left of the particle
with label $J(t)$. Hence $G(J(t), I(t))\leq t$, but $G(J(t)+1, I(t))>t$
and $G(J(t), I(t)+1)>t$. The next jump in the process $(I(t), J(t))$
occurs either when particle with label $J(t)$ jumps right, or when the 
hole with label $I(t)$ moves left, whichever happens sooner. 
Let $\phi_n$ be the $n$th term in the jump chain of the process $(I(t), J(t))$. 
Then $\phi_0=(1,1)$ and we have the following recurrence for the sequence $\phi_n$:
\begin{equation}
\label{phirec}
\phi_{n+1}=
\begin{cases}
\phi_n+(0,1)&\text{if } G\big(\phi_n+(0,1)\big)<G\big(\phi_n+(1,0)\big)
\\
\phi_n+(1,0)&\text{if } G\big(\phi_n+(1,0)\big)<G\big(\phi_n+(0,1)\big).
\end{cases}
\end{equation}
So in terms of the growth model, $\phi_n$ behaves as follows. We start from 
the initial condition in which only site $(1,1)$ is occupied, and $\phi_0=(1,1)$.
Then at each step, $\phi_n$ moves from its current position to whichever of its
neighbours above and to the right is the first to be occupied.

We now describe how the sequence $\phi_n$ corresponds to a \textit{competition interface}
in the same growth model when occupied sites take one of two colours.
The rules are:
\begin{itemize}
\item Any site $(x,1)$, $x>1$ takes the colour blue.
\item Any site $(1,x)$, $x<1$ takes the colour red.
\item For $x,y>1$, the site $(x,y)$ takes the same colour as 
whichever of its neighbours $(x-1,y)$ and $(x,y-1)$ became occupied
most recently.
\end{itemize}
The colour of $(1,1)$ is not important. 

Let $\cB$ be the set of sites that take the colour blue (when they are eventually occupied)
and let $\cR$ be the set of set of sites that take the colour red. 
So $\cR\cup\cB=\NN^2\setminus\{(1,1)\}$. 
Using a simple induction argument based on the growth rules above and on (\ref{phirec}),
we get the following characterisation
which justifies the name ``competition interface'': 
the set $\{\phi_n, n\geq0\}$ consists precisely of those points 
$z$ such that $z+(0,1)\in\cR$ while $z+(1,0)\in\cB$.

In addition, if $(x,y)$ is in the competition interface (i.e.\ if 
$(x,y)=\phi_n$ for some $n$, or equivalently if $(x,y)=\big(I(t), J(t)\big)$ for 
some $t$)
then all the points $\{(x+k, y), k>0\}$ to its right are in $\cB$ and 
all the points $\{(x, y+k), k>0\}$ above it are in $\cR$.

\begin{figure}[h]
\centering
\includegraphics[width=0.5\linewidth]{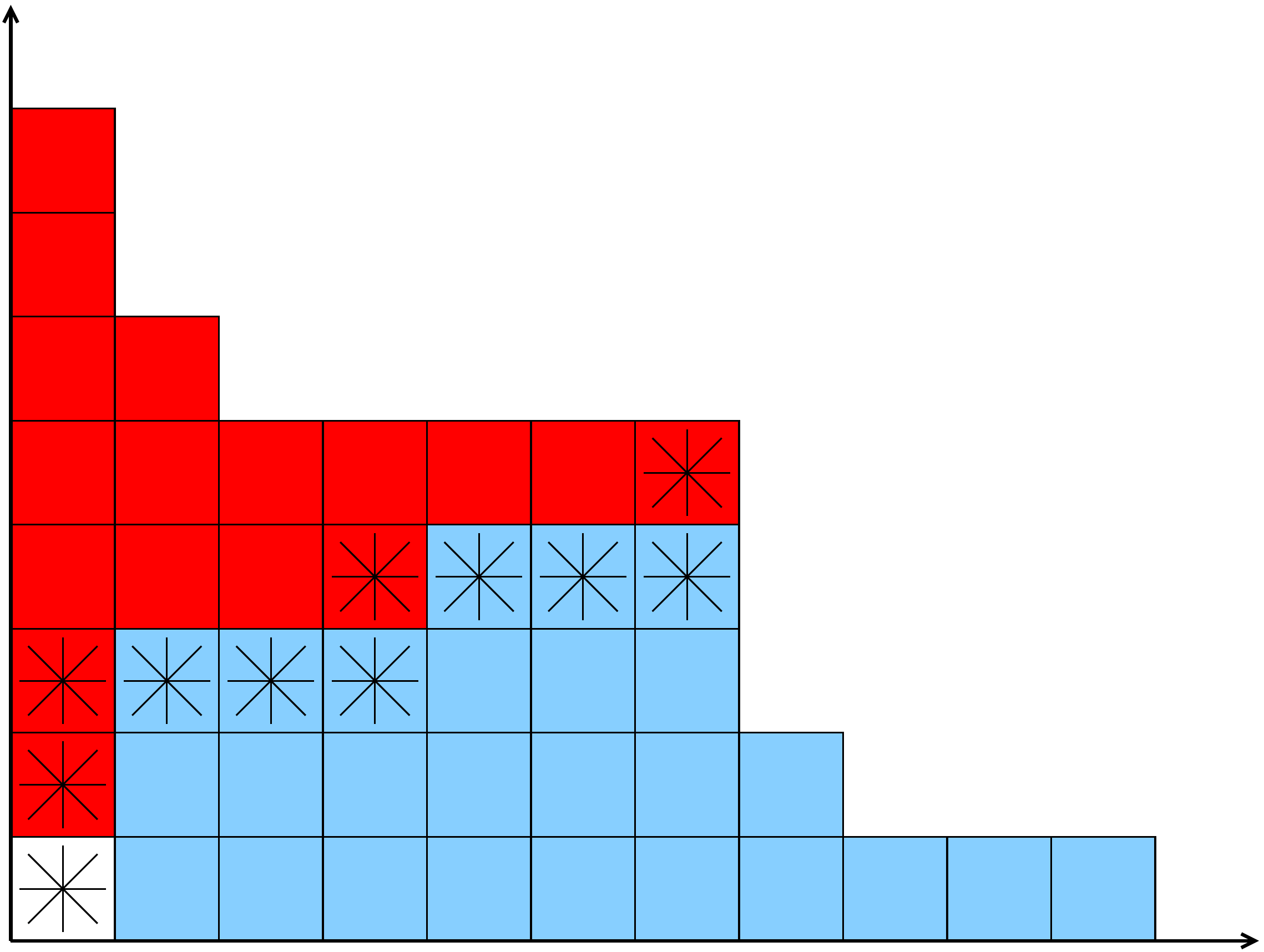}
\caption{\label{oneinterfacefig}
An example of the evolution of the growth model with one competition interface.
The asterisks indicate those already occupied sites which form
part of the competition interface. For example, since the 
competition interface passes from $(1,1)$ to $(1,2)$, $(1,3)$ and $(2,3)$, 
we see that $G(1,2)<G(2,1)$, that $G(1,3)<G(2,2)$ and that $G(2,3)<G(1,4)$.}
\end{figure}

See Figure \ref{oneinterfacefig} for an example of the competition interface.

Now we are finally ready to turn to the \TASEP with \textit{two} second-class
particles, and the corresponding coexistence question in Theorem \ref{growththm}.
In the same way that the system with one second-class particle corresponded to a growth
model with one interface, the system with two second-class particles will correspond
to the growth model with two interfaces.

Consider the \TASEP starting from the state
$\hat\parti\un\second\second\hat\hole$.  We can again replace the sites
containing second-class particles by pairs of sites, each containing a hole and
a particle in that order.  This representation works just as it did in the
system with a single second-class particle, \textit{until the moment when one
  second-class particle tries to jump onto the other}.  At this point we
  move from $[\hole\parti][\hole\parti]$ to $\hole\hole\parti\parti$ and we can
  no longer maintain the two $[\hole\parti]$ pairs (although in fact 
one can identify a
  *pair in the configuration $\hole[\hole\parti]\parti$ and show that this *pair
  is in correspondence with a competition interface -- see the remark below).  
In any case, up until the
collision moment, the two second-class particle representation is valid, and
in particular Theorem \ref{t1}(a) tells us that the probability that such a 
collision is ever attempted is $2/3$.

\newcommand{\Ileft}{\overleftarrow{I}}
\newcommand{\Jleft}{\overleftarrow{J}}
\newcommand{\Iright}{\overrightarrow{I}}
\newcommand{\Jright}{\overrightarrow{J}}
\newcommand{\phileft}{\overleftarrow{\phi}}
\newcommand{\phiright}{\overrightarrow{\phi}}
\newcommand{\cBleft}{\overleftarrow{\cB}}
\newcommand{\cBright}{\overrightarrow{\cB}}

We will again label the particles and holes sequentially.
The initial state is:
\begin{center}
\includegraphics[width=0.6\linewidth]{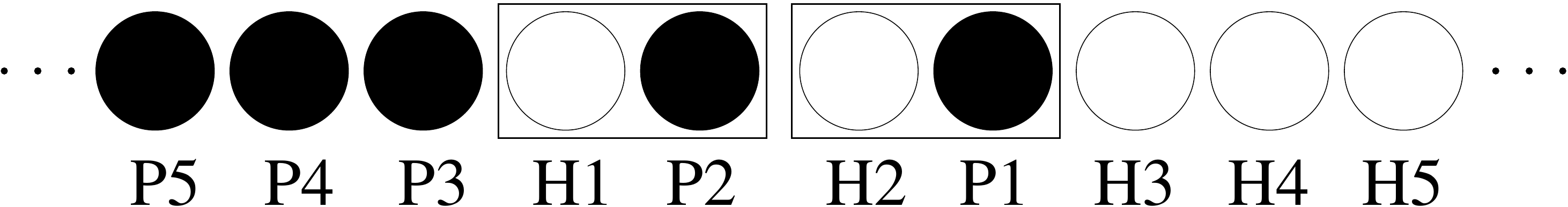}
\end{center}
In this state particle P1 has overtaken holes H1 and H2, and particle P2 has
overtaken hole H1. Hence we put $G(1,1)=G(2,1)=G(1,2)=0$.

We will now keep track of both the \starpairs:
let $\Ileft(t)$ and $\Jleft(t)$ be respectively the labels of 
the hole and the particle covered by the left-hand \starpair,
and let $\Iright(t)$ and $\Jright(t)$ be those covered
by the right-hand \starpair. 
We have $\big(\Ileft(0), \Jleft(0)\big)=(1,1)$
and $\big(\Iright(0), \Jright(0)\big)=(2,1)$.
More generally, if at time $t$ 
the two \starpairs
are next to each other at time $t$ if and only if 
\begin{equation}\label{nexteq}
\Iright(t)=\Ileft(t)+1 \text{ and } \Jleft(t)=\Jright(t)+1.
\end{equation}
The local evolution both of $\left(\Ileft(t), \Jleft(t)\right)$ 
and of $\left(\Iright(t), \Jright(t)\right)$ depends on 
the values $G(z)$ in just the same way as in the paragraph 
around (\ref{phirec}).
If we define $\left(\phileft_n, n\geq0\right)$ and
$\left(\phiright_n, n\geq0\right)$ to be the jump chains
of these two processes, then each of them separately obeys
the same recurrence as the process $\left(\phi_n\right)$ at (\ref{phirec}).

Now recall the rules of the growth model with two interfaces from Section
\ref{growthsec}. In the initial condition,
the sites $(1,1)$, $(1,2)$ and $(2,1)$, and all other sites are unoccupied. 
When a new site becomes occupied, its colour is chosen as follows:
\begin{itemize}
\item Any site $(1,x)$ becomes dark blue.
\item Any site $(x,1)$ becomes light blue.
\item The site $(2,2)$ becomes red.
\item A site $(x,y)$, $x>1$, $y>1$, $(x,y)\ne(2,2)$, takes the colour of whichever
of its parents became occupied \textit{more recently}.
\end{itemize}
\paragraph{\sl Remark} In fact, we have elaborated 
slightly by introducing both dark and light blue
colours at this point for convenience, to distinguish the clusters on either
side of the red cluster. The difference is convenient for 
the argument below, but it is not significant for the
result, since we only ask whether or not the red cluster is unbounded, and we
are free to identify the two blue clusters as before whenever we
like. In addition, once we make the distinction between the two blue clusters, 
we can consider the interface between those clusters \emph{after} the red
cluster disappears. This interface is identified
by the *pair resulting from the coalescence of the *pairs related to the
light-blue/red and red/dark-blue interfaces,
if we make the jump from $[\hole\parti][\hole\parti]$ 
to $\hole[\hole\parti]\parti$ suggested above.

As before, the processes $\left(\phileft_n\right)$ and $\left(\phiright_n\right)$
are competition interfaces, with the following properties. 
Let $k\geq1$ and $n\geq0$. Then any point $\phiright_n+(k,0)$ to the right of the 
lower/right interface is (light) blue, and any point $\phiright_n+(0,k)$ 
which is above that interface is either red or (dark) blue.
Similarly any point $\phileft_n+(0,k)$ which is above the upper/left
interface is (dark) blue, and any point $\phileft_n+(k,0)$ which is to the 
right of that interface is either red or (light) blue.

\begin{figure}[h]
\centering
\includegraphics[width=0.5\linewidth]{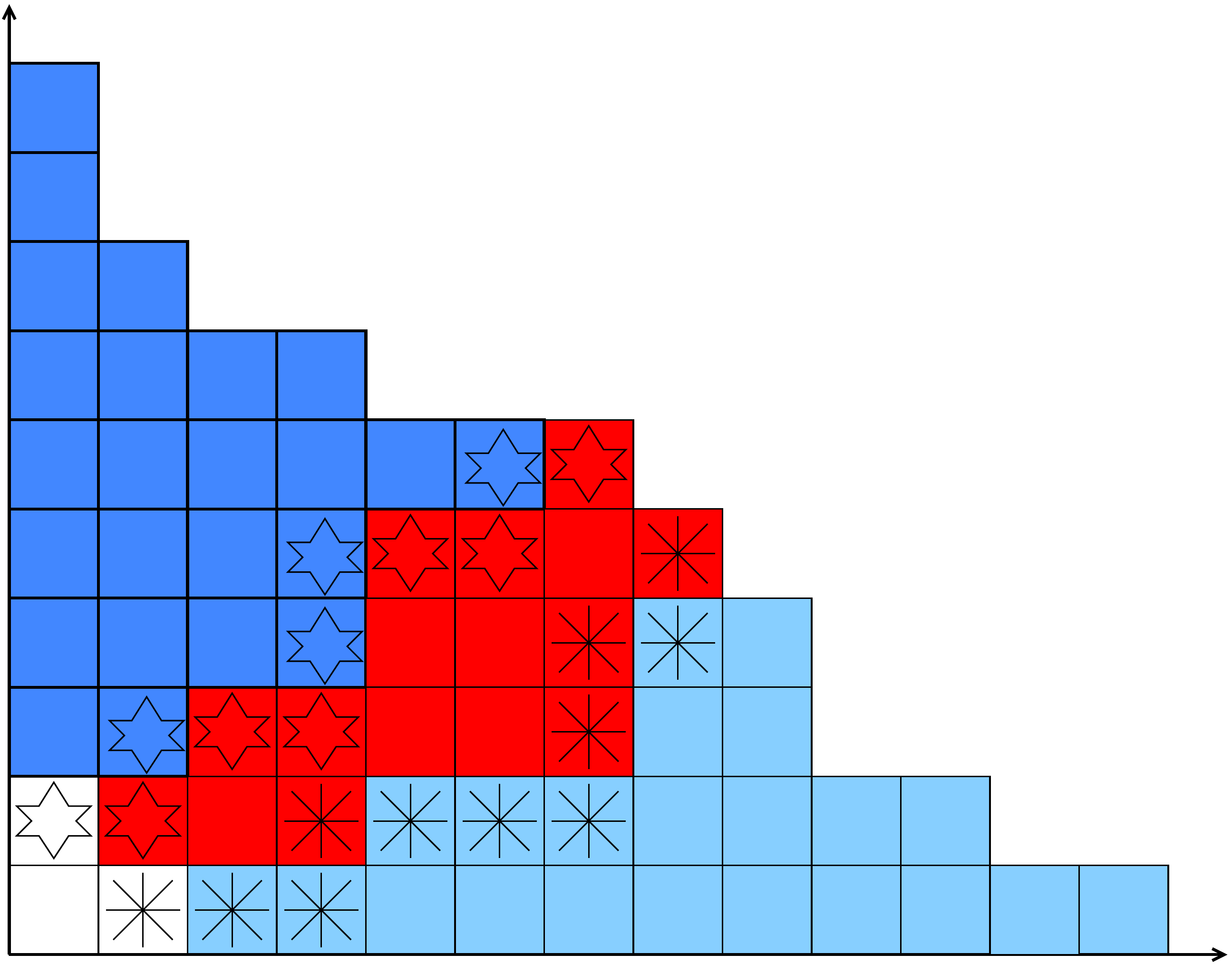}
\caption{\label{meetingfig} An example of the evolution of the
  growth model with two interfaces.  Here we see that
  $G(3,1)<G(2,2)<G(1,3)$, since the lower-right interface (indicated by
  asterisks) moved to $(3,1)$ rather than $(2,2)$ but the upper-left interface
  (indicated by stars) moved to $(2,2)$ rather than $(1,3)$.  In the state
  shown, the two interfaces are at $(8,5)$ and $(7,6)$.  If $G(8,6)$ is smaller
  than both $G(9,5)$ and $G(7,7)$ then both interfaces will move next to the
  point $(8,6)$.  In this case the point $(8,6)$ will be red but no further
  points will ever become red; the red cluster will be surrounded.  If, instead,
  one of $(9,5)$ or $(7,7)$ is occupied before $(8,6)$, the red cluster will
  grow further, perhaps unboundedly.}
\end{figure}

From (\ref{nexteq}), we see that 
the moment (if it ever happens) when the two interfaces coincide for the first time
is the first moment when the left-hand second-class particle attempts to jump 
onto the right-hand second-class particle. 
Once the interfaces meet in this way, the red cluster is ``surrounded'' and cannot
grow any further. Specifically, the last paragraph tells us that all the points above
the upper-left interface will be (dark) blue and all those to the right of the 
lower-right interface will be (light) blue, so that no further points will be able to take on the 
colour red. See Figure \ref{meetingfig} for an example. 

On the other hand, suppose that the second-class particles have not yet collided, 
so that the interfaces have not yet coincided. Even if the interfaces are already 
adjacent, it is still possible for them to move apart again, and any sites
between the two interfaces must take the colour red; so certainly the 
red cluster is not yet surrounded. 
See Figure \ref{meetingfig} for an illustrative example.

Hence the event that the red cluster grows unboundedly corresponds to the
event that the two second-class particles never meet. From Theorem
\ref{t1}(a), we know that that event has probability $1/3$.

\section*{Acknowledgments} We thank Soledad Torres and Karine Bertin for
nice discussions about their simulation results. We are most grateful to a 
referee who pointed out an important error in an earlier version of the paper. 

P.A.F.\ is partially supported by CNPq and FAPESP grants. 
P.G.\ was supported by FAPESP--Brasil with the grant 06/58527-0.

\end{document}